\documentclass[11pt]{amsart}
\usepackage{latexsym}
\usepackage{amsfonts,amsmath,amssymb}

\usepackage{graphicx}
\usepackage{amsfonts}
\usepackage{hyperref} 
\usepackage{textcomp}
\usepackage{amssymb}

\newtheorem{theorem}{Theorem}
\newtheorem{definition}[theorem]{Definition}
\newtheorem{lemma}[theorem]{Lemma}
\newtheorem{conjecture}[theorem]{Conjecture}
\newtheorem{proposition}[theorem]{Proposition}
\newtheorem{corollary}[theorem]{Corollary}

\newtheorem{question}{Question}

\newtheorem{example}{Example}

\begin{document}
\author{ Mikl\'os B\'ona 
         \and Ryan Flynn}
\thanks
{\noindent Department of Mathematics\\
University of Florida\\
Gainesville FL 32611-8105}

\title{Sorting a Permutation with Block Moves}

\begin{abstract}
We prove a lower and an upper bound on the number of block moves necessary
to sort a permutation. We put our results in contrast with existing results
on sorting by block transpositions, and raise some open questions.
\end{abstract}

\maketitle
\section{Introduction}

\subsection{Background and Definitions}
Let $p=p_1p_2\cdots p_n$ be a permutation. A {\em block} of $p$ is just a
string of consecutive entries $p_ip_{i+1}\cdots p_j$. 
A {\em block transposition}
is the operation that interchanges two consecutive blocks of $p$, without
changing the order of entries within each block. So if the two blocks
are $p_ip_{i+1} \cdots p_{j-1}$ and $p_{j}p_{j+1}\cdots p_{k-1}$, then the
block transposition interchanging these two blocks results in the
permutation
\[p_1p_2\cdots p_{i-1}p_jp_{j+1}\cdots p_{k-1}p_ip_{i+1}\cdots p_{j-1}
p_k\cdots p_n.\]

There are some interesting results and intriguing conjectures on sorting
permutations with block transpositions. The interested reader should consult
\cite{wastlund} for an overview.  The original motivation came from 
molecular biology and genome sorting. Details and examples of how these moves
 occur during evolutionary processes can be found in \cite{bafna}.

In this paper, we consider a similar, but not identical class of operations,
called {\em block moves}. A block move still interchanges two blocks of
entries, but the two blocks do not have to be adjacent anymore. It can
happen that one block is $p_{i}p_{i+1}\cdots p_{i+a}$, while the other block
is $p_{j}p_{j+1}\cdots p_{j+b}$, where $i+a<j-1$. 
In particular, every block transposition is a block move, but not every 
block move is a transposition.

We will show that while in general, sorting by block moves seems to be 
much more efficient than sorting by block transpositions, for some permutations
that are very hard to sort, there is only a very small difference between
the efficiency of the two sorting algorithms. 

We also show that every permutation of length $n$ (or, in what follows, 
$n$-permutation) can be sorted by at most $\lfloor (n+1)/2 \rfloor$
 block moves. The corresponding statement is only conjectured for sorting
by block transpositions \cite{wastlund}.

\subsection{Earlier results}

If a series of operations takes a permutation $p$ 
to the increasing permutation $12\cdots n$, then we say that that series
of operation {\em sorts} $p$. It is natural to ask how many operations of a
given kind are necessary to sort a given permutation $p$. For block
transpositions, the best result is given in \cite{wastlund}.

\begin{theorem} \cite{wastlund} 
Let $n\geq 9$. Then every $n$-permutation can be sorted by at most
$\lfloor (2n-2)/3 \rfloor $ block transpositions.
\end{theorem}

On the other hand, it is proved in \cite{wastlund} that to sort the
decreasing permutation $n(n-1)\cdots 21$, one needs exactly
$\lceil (n+1)/2 \rceil $ block transpositions, for all $n\geq 3$. This leaves
the intriguing question as to where between $n/2$ and $2n/3$ the actual
number of needed block transpositions lies. In response to this question,
the authors of \cite{wastlund} stated the following conjecture.

\begin{conjecture} If $n\neq 13$ and $n\neq 15$, then every $n$-permutation
can be sorted by at most $\lceil (n+1)/2 \rceil $ block transpositions.
\end{conjecture}
If this conjecture is true, that  means that no permutation is harder
to sort by block transpositions than the decreasing permutation. 

In this paper, we attack the analogous questions for block moves instead
of block transpositions.

\section{Sorting by Block Moves}
\subsection{A Lower Bound}
A {\em descent} of a permutation $p=p_1p_2\cdots p_n$ is an index $i$ so
that $p_i>p_{i+1}$. For instance, $p=34152$ has two descents, $i=2$ and 
$i=4$. Our main tool in this subsection is the following lemma.

\begin{lemma} No block move decreases the number of descents of a permutation
$p$ by more than two.
\end{lemma}

\begin{proof}
Let us assume the contrary. That is, let us assume that there exists a 
permutation \[p=p_1\cdots s\underline{x\cdots y} u \cdots v 
\underline{w\cdots z}t\cdots p_n\] so that upon
 interchanging the two underlined
blocks $\underline{x\cdots y}$ and $\underline{w\cdots z}$, the obtained 
permutation \[p'=p_1\cdots s \underline{w\cdots z}  u \cdots v 
\underline{x\cdots y} t\cdots p_n\] has at least three less descents 
than $p$ does.

Clearly, if $i$ is a descent of $p$, then $i$ is a descent of $p'$ unless
$p_i$ and $p_{i+1}$ are separated by the block move that turns $p$ into $p'$.
That is, the descents that could possibly be destroyed in the displayed 
example are at the positions of $s$, $y$, $v$, and $z$ in $p$. This is four 
positions. If descents in at least three of them are destroyed, then there
 are two ways in which  that
can happen. Either all four were descents in $p$, but only at most one of those
positions is a descent in $p'$, or   three of them were
 descents in $p$, and none are
descents in $p'$. Without loss of generality, we can assume that we are
in the first case.

 Indeed, if we are in the second case, then $p$ has
three descents in the four considered positions, and $p'$ has  none.
Let $p^r$ be
the reverse of $p$, and let $p'^r$ be the reverse of $p'$. Taking reverses
turns non-descents into descents and vice versa. So $p'^r$ now has four
descents
in the considered positions, and $p^r$ has one. As $p'^r$ can certainly be 
taken to $p^r$ by one block move, this reduces the second case to the
first.

So assume that $s>x$, $y>u$, $v>w$, and $z>t$, and also, in $p'$, there is no
descent in at least three of these four positions.
\begin{enumerate}
\item  Let us say these are the
 first three positions. That is, assume that $s<w$, $z<u$, and $v<x$. This is
a contradiction, since together with the inequalities assumed in $p$,
we get the chain of inequalities
$s<w<v<x<s$.
\item Now assume that descents are destroyed in the first, third, and fourth
considered positions. That means $s<w$, $v<x$, and $y<t$. 
This again leads to the chain of inequalities  $s<w<v<x<s$.
\item If descents are destroyed in the first, second, and fourth considered 
positions, then
 $s<w$, $z<u$, and  $y<t$. Together with the inequalities assumed
in $p$, this leads to $t<z<u<y<t$.
\item Finally, if descents are destroyed at the second, third, and fourth 
positions, we have $z<u$, $v<x$, and $y<t$. This again leads to 
the chain of inequalities $t<z<u<y<t$.
\end{enumerate}
This shows that the number of descents cannot decrease by more than two 
during any block move. 
\end{proof}

\begin{corollary} \label{lowerbound}
If a permutation $p$ has $d(p)$ descents, then at least $\lceil d(p) /2
\rceil $ block moves are needed to sort $p$. 
\end{corollary}

In particular, $\lceil (n-1)/2 \rceil $ block moves are needed to sort the
decreasing $n$-permutation $n(n-1)\cdots 21$. It is easy to see that that 
number of moves is actually sufficient. One can just interchange
the one-element block $i$ and the one-element block $n+1-i$, for all 
$i\leq \lceil (n-1)/2 \rceil $. 

It is proved in \cite{wastlund} that one needs exactly $\lceil (n+1)/2 
\rceil $ block {\em transpositions} to sort the decreasing $n$-permutation.
Our Corollary \ref{lowerbound} shows that though block moves are much more
general than block transpositions, in the worst case (that of decreasing
permutations), they are barely more efficient.

\subsection{An Upper Bound}

In the previous subsection, we showed what sorting by block moves {\em
cannot} accomplish; now we will show what it {\em can}. 
That is, we will show that it can sort all $n$-permutations in
roughly $n/2$ steps. 

In the previous subsection, the number of descents of a permutation $p$
was used to find a lower bound for the number of block moves needed to sort
$p$. It seems that the number of descents is not the right statistic
to use when looking for more refined results. For instance,
permutations 21345678 and 24681357 both have one descent, but it is clear
that the former is much easier to sort than the latter. 

Therefore, we will resort to another statistic to measure the un-sortedness of
a permutation.

\begin{definition}
Let $p=p_1p_2\cdots p_n$ be a permutation. 
We say that the pair $(i,i+1)$ of indices (with $0\leq i\leq n$)
 is a \textit{good pair}
 in  $p$ if any of the following conditions hold:
\begin{enumerate} \item $1\leq i\leq n-1$ and $p_i+1=p_{i+1}$, or
\item $i=0$ and $p_1=1$, or
\item $i=n$ and $p_n=n$.
\end{enumerate}

 A \textit{bad pair} is any pair of consecutive  indices $(j,j+1)$ with
$0\leq j\leq n$ that
 is not a good pair.
\end{definition}

 For instance, 
in $41253$, the indices 2 and 3 form a good pair, and there are no other
good pairs in this permutation. Note that if $p$ is not the increasing
permutation, then $p$ always has at least two bad pairs. 

One measure of the un-sortedness of a
 permutation is the number of bad pairs. We will use this number to prove
an upper bound on the number of block moves needed to sort a particular 
permutation.

For  $p=p_1p_2\cdots p_n$, we will  sometimes write 
 $p=0p_1p_2\cdots p_n(n+1)$. This simply illustrates the fact
 that moving the entry 1 to the beginning of $p$ or the entry $n$
 to the end of $p$ actually decreases the number of bad pairs by one.

\begin{proposition}
If $p$ is not the increasing permutation, 
 then there exists a block move that decreases the number of bad pairs of
$p$ by at least two. 
\end{proposition}

\begin{proof}
First, note that we can assume that $p$ has no good pairs. Indeed, 
  if it did, then one could glue the two entries forming that good pair
 together and consider the obtained $(n-1)$-permutation instead.
  Even after this reduction, the length of $p$ is at least two,
 since otherwise $p$ would have been the increasing permutation.

Similarly, we can assume that $p_1 \neq 1$. We can also assume that 
$p_1 \neq 2$, since otherwise  there would be a move reducing the number
of bad pairs by two - namely the move switching $p_i=1$ and 
the block $p_1p_2\cdots p_{i-1}$.

So, without loss of generality, we can assume that  $p_1 \geq 3$.  

In a similar manner, we can assume that $p_n\neq 1$, since otherwise
 the move interchanging  the block consisting of the entry $p_n=1$
 with the block  $p_1p_2\cdots p_j$, where $p_{j+1}= 2$
 would remove two bad pairs. We can proceed analogously in the more general
case when the entry 2 precedes the entry 1 in $p$. 

Consider all decreasing subsequences of $p$ starting at $p_1$; that is,
 all sequences
$p_1 > p_{i_1} > p_{i_{2}} \cdots > p_{i_k}$, 
where $1 < i_{1} < i_{2} < \cdots < i_{k}$.

There is at least one such sequence since $p_1 \geq 3$.

Among these sequences, consider those of the longest horizontal length;
 that is, those which terminate as far to the right as possible. If there are
several candidates, choose the one in which no entry can be replaced by an
an entry on its left (and keeping all other entries unchanged), and which
 cannot be extended by inserting an additional entry. Call this subsequence
$p_{dec}$. In general, if a decreasing sequence $seq$
 has the property that none
of its entries can be replaced by an entry on its left so that if all other 
entries are left unchanged, the obtained sequence is still decreasing, then
we will say that $seq$ has no {\em left refinement}.
 
\begin{example}
Let $p=51342$. Then there are two decreasing subsequences of $p$ that start
at $p_1=5$ and end at $p_5=2$, namely 542 and 532. We choose 
$p_{dec}=532$, since in 542,
the second entry, 4, can be replaced by one on its left, 3. So 542 does have
a left refinement, but 532 does not. 
\end{example}

 Now note that each entry on the right of  $p_{i_{k}}$ is larger than all
 the $p_{i_j}$, since otherwise  we could increase the horizontal length of
 our decreasing sequence $p_{dec}=
p_1p_{i_1}\cdots p_{i_k}$.  So, the entry 1 must
 be in one of the blocks  $p_{i_j}\cdots  p_{i_{j+1}}$.  
(If $1=p_{i_k}$,then it is easy to show that 2 preceeds 1 in $p$,
 since $p_1 \geq 3$.)
 
 If $p_{i_j} - 1 = p_{i_{j+1}}$, then we can find a block
 move that decreases the
number of bad pairs by two, namely the block move switching the two
 underlined blocks below:
 
 \[p = 0 \underline{p_1p_2\cdots  p_{i_{j}-1}} p_{i_{j}}\cdots
 \underline{1 ... p_{i_{j+1}}} \cdots .\]

 So, we may assume that $p_{i_{j+1}} \leq p_{i_{j}} - 2$.
 
 Now consider $p_{i_{j}} - 1$.  If this entry precedes $p_{i_{j+1}}$, then
 it could replace it in $p_{dec}$. So, $p_{i_{j}} - 1$ must lie to the right 
of $p_{i_{j+1}}$.  Then  the following block move decreases the number of
 bad pairs by two:
 
\[ 0 \underline{ ..... p_{i_{j} - 1}} p_{i_{j}} \cdots
 \underline{1 . . . p_{i_{j+1}}\cdots p_{i_{j}} - 1} \cdots .\]
\end{proof}

\begin{example}
Let $p=  0 \  \underline{9}\ 3\ 10\ 6 \ \underline{8}\ 2\ 4\ 1\ 5 
\ \underline{7} \ 12 \ 11  \  13$.
Then 9 8 7 is a decreasing sequence of maximal horizontal length, and has no
left refinement. The entry 1 is in the block between 8 and 7.  In this 
case, 8-1 = 7, so by the proof above the block move switching the two 
underlined blocks below
\[0 \ \underline{9 \ 3 \ 10 \ 6}\ 8 \ 2 \ 4 \ 
\underline{1\ 5\ 7}\ 12\ 11 \  13\]
decreases the number of bad pairs by two.
\end{example}

The main result of this section now follows.

\begin{theorem} Let $p$ be a permutation, and let $b(p)$ be the number
of bad pairs of $p$. Then $p$ can be sorted by $b(p)/2$ block moves. 

In particular, every $n$-permutation can be sorted by at most $\lfloor
(n+1)/2\rfloor$ moves. 
\end{theorem}

\begin{proof} Start sorting $p$ and in each step, use a block move that
decreases the number of bad pairs by at least two. The preceding lemma shows
that there is always such a move.
\end{proof}

\section{Further Directions}

Unlike the number of descents, the number of bad pairs {\em can} go
down by three or four in one block move. For instance, if $p=14325$, then
$p$ has four baid pairs, but interchanging the second entry and the fourth
entry results in the identity permutation, that has zero bad pairs. 

The question is how often will this happen. It is not difficult to see
that in a random permutation, there is a high probability that there is a
block move reducing the number of bad pairs by more than two. However, it
is less clear how many such moves we will find during the process of 
sorting a permutation, and whether it pays to be greedy, that is, to
always choose the block move that reduces the number of bad pairs by the
largest number. 

These observations suggest the following questions.

\begin{question} \label{ambitious} 
Is there a permutation statistic $z(p)$ that is easy
to define directly and that allows for direct computation of the number
of block moves necessary to sort $p$?
\end{question}

Our results in this paper show that at the very least, $\max (d(p)/2, b(p)/4)$
block moves are needed to sort $p$. 

\begin{question}
Is there a positive constant $c<0.5$ so that if $n$ is sufficiently large,
then almost all $n$-permutations can be sorted by no more than $cn$ block
moves?
\end{question}

\begin{question}
Does a greedy approach work? That is, if we always choose a move that 
decreases the number of bad pairs by as much as possible, will we sort
our permutation as fast as possible?
\end{question}

Numerical evidence suggests that there are many $n$-permutations that are
as difficult as possible, that is, that take $\lfloor (n+1)/2 \rfloor $
block moves to start. This raises the following question that may be easier
than the similar Question \ref{ambitious}.

\begin{question} Which $n$-permutations take the maximal 
 $\lfloor (n+1)/2 \rfloor $ block moves to sort? How many such permutations
are there?
\end{question}

Finally, we point out that the distribution of the number of good pairs in 
a random $n$-permutation is fairly well-understood. It has been proved in
\cite{wolfowitz} that this distribution converges (in distribution)
 to a Poisson distribution
of parameter 1 as $n$ goes to infinity. See \cite{bona} for an overview. 
So most $n$-permutations will have only a few good pairs in them. So 
improvements on the number of necessary moves will not come from taking a
better look at the starting permutations; they would have to come from 
taking a more detailed look at the sorting methods.


\begin{thebibliography}{99}
\bibitem{bafna} V. Bafna, P. Pevzner, Sorting by Transpositions, {\em 
SIAM J. Discrete Math.} {\bf 11} (1998) 224--240.
\bibitem{bona} M. B\'ona, On Three Notions of Monotone Subsequences, 
{\em Proceedings of the London Mathematical Society}, to appear. 
\bibitem{wastlund} H. Eriksson; K. Eriksson; J. Karlander; L. Svensson; J.
 Wästlund, Sorting a bridge hand. Selected papers in honor of Helge 
Tverberg.  {\em Discrete Math.}  {\bf 241}  (2001),  no. 1-3, 289--300. 
\bibitem{wolfowitz} J. Wolfowitz, Note on Runs of Consecutive Elements,
{\em Annals Math. Statistics,} {\bf 15} (1944), 97-98.
\end{thebibliography}
\end{document}